\documentclass[12pt,leqno,a4paper] {amsart}
\usepackage{amssymb,enumerate}
\usepackage{hyperref}
\usepackage{csquotes}

\newcommand{\CC}{{\mathbb{C}}}

\newcommand{\ZZ}{{\mathbb{Z}}}
\newcommand{\cO}{{\mathcal{O}}}

\newcommand{\fS}{\mathfrak{S}}
\newcommand{\fg}{\mathfrak{g}}

\theoremstyle{definition}

\begin{document}

\title[Carter]{Roger Carter}

\date{\today}
\author{Meinolf Geck}
\address{Universität Stuttgart,
Fachbereich Mathematik,
IDSR - Lehrstuhl für Algebra,
Pfaffenwaldring 57,
D-70569 Stuttgart, Germany.}
\email{meinolf.geck@mathematik.uni-stuttgart.de}
\author{Donna M. Testerman}
\address{Institut de Math\'ematiques, Station 8, \'Ecole Polytechnique
  F\'ed\'erale de Lausanne, CH-1015 Lausanne, Switzerland.}
\email{donna.testerman@epfl.ch}


\thanks{}




\maketitle


Roger Carter (1934--2022) was a very well known mathematician working in
algebra, representation theory and Lie theory. He spent most of his mathematical career in Warwick. Roger was
a great communicator of mathematics: the clarity, precision and enthusiasm
of his lectures delivered in his beautiful handwriting were hallmark
features recalled by numerous students and colleagues. His books
have been described as marvelous pieces of scholarship and service to the
general mathematical community. We both met Roger early in our careers, and were encouraged and influenced
by him~---~and his lovely sense of humour. This text is our tribute, both
to his mathematical achievements, and to his kindness and generosity
towards his students, his colleagues, his collaborators, and his family.

\section{His Life}

\centerline{\it{Early life}\footnote{Many thanks to Audrey Edwards for sharing with us  a text which she wrote after Roger Carter's death, much of the content of which is incorporated here.}}\bigbreak

\noindent Roger Carter was born 25 August 1934, in Shipley, West Yorkshire. His father was employed by the Inland Revenue as an Income Tax Inspector
 and because of this employment, the family moved often.  The first home for Roger and
his sister Audrey (his junior by 17 months), was Baildon, Yorkshire and the second home was Workington, in Cumberland. Roger's father always
insisted that the frequent  moves should not interrupt Roger's schooling, so his junior school was in Chesterfield and after he passed the exam admitting
him to grammar school the
family moved to
Barnsley where Roger attended the Barnsley Boys' Grammar School. He excelled in mathematics and was duly awarded a scholarship to Sidney Sussex College, Cambridge.

Roger's sister relates that Roger always had a sense of humour and the makings of the excellent lecturer he would become: when he was learning algebra
at school, he decided to teach his younger sister. He said to her ``Write this down - Let $x$ be a pair of large red knickers and let $y$ be a
pair of small
blue knickers...'' His sister recalls that the lesson didn't progress very far because she was giggling too much. Years later, Roger showed the same pedagogical sense of humour when participating in Miranda Mowbray's PhD exam:
\begin{displayquote}{\small ``Roger Carter was the external examiner for my PhD. At one point in
the oral exam he asked me ``can you think of a Lemma that could be applied at this part of your thesis?'' I felt nervous, thought for a while and then said
``I don't know, I can't think of anything''. He leaned forward and said, ``Are you SURE?'' and I said ``Oh! You could apply Schur's Lemma.'' After this, I felt a lot less nervous.''}
\end{displayquote}

In addition to arranging that Roger's schooling progressed smoothly, his father made sure that his children became committed and sturdy walkers.
He sold his car after their births and indeed Roger and Audrey both came to thoroughly enjoy walking. Roger's colleagues will recall his long
stride on the Yorkshire
Moors when the excursion day at a Durham meeting was organised. George Lusztig, Roger's former colleague and collaborator, remembers hiking with Roger in the Zillertal Alps in Austria,
around 1973,
while also discussing mathematics.

During the war years, the family lived in Chesterfield, and thereby escaped the bombing, and were only sometimes woken by planes flying over to bomb the
Sheffield Steelworks. Roger and his sister played together a lot during that time, with Roger mostly taking the lead, as Audrey describes him as
being ``much braver''.
Then when they were
older, they continued playing together, also as a family, board games and card games. Even in retirement, Roger took his monopoly games seriously.

Roger's father also taught his children to row and to play chess, and Roger pursued both activities as an adult. He rowed for his college in
Cambridge and became an excellent chess player.
Roger was interested in other sports as well, enjoyed refereeing football and attended football matches and cricket games, sometimes taking his
sister along. Another favorite leisure activity
in which
Roger partook was
singing. He had a lifelong love for classical music, sang in the church choir wherever he lived, and while in Cambridge sang with the
Cambridge University Musical Society. He also
played the piano and for some years was the organist at the Holy Trinity Church in Leamington Spa. This love of music and his general
good humour and playfulness added a special note to the 1987 Durham Mathematical Symposium
which he co-organised. After one of the
evening meals,
Roger and a small set of conference participants prepared a choral concert, and led the full group of participants  in singing two
``mathematical songs'', for which
they had written ad-hoc lyrics based upon the mathematical topic of the
conference and its participants' names,
sung to the tunes of ``John Brown's Body'' and ``Do-Re-Mi''.

After the war, the family regularly took holidays in Scarborough. As teenagers, Roger and his sister spent a number of Christmas holidays at a
Quaker guest house north of Scarborough and were
later offered a free holiday in exchange for organising all of the entertainment for the 100 guests. Roger was in charge of table tennis and
bar billiards
competitions and his sister set up dancing games. And again, there was a concert. Learning of this, one is then not surprised that Roger took
pleasure in organising international conferences, research programmes and Warwick-based thematic symposia, and in hosting visitors to Warwick. Roger often spent his
Christmas holidays with his sister and, after she moved to
Scarborough, he visited her family for the annual cricket festival as well, and would spend days out with a map exploring the coast path. His niece Helen Smith
says that her chief memories are of ``someone who arrived punctually each year for Christmas, very much enjoyed teasing his young great nephews and nieces and
had a tremendously sweet tooth and a passion for licorice allsorts! Really, mathematics was his life.''
\bigbreak

\centerline{\it{Cambridge, T\"{u}bingen, Newcastle, Warwick}}\bigbreak

\noindent Roger matriculated at Cambridge University in 1952, where he obtained a first class degree in mathematics.
His sister says that the arrival in Cambridge for his undergraduate degree was a challenge, as he found himself a Yorkshire
``grammar school boy'' among crowds of public school boys on the railway platform. Nevertheless, he settled in and thrived academically and socially, and
went on to write a
PhD thesis under the direction of Derek Taunt. The title of his thesis, defended in 1959, was ``Some Contributions to the Theory of Finite Soluble
Groups''. Roger was elected a member of the London Mathematical Society 18 June 1959 and remained a member throughout his life. 

After obtaining his PhD, he spent a year at the Univerity of T\"{u}bingen on a Humboldt Stiftung research fellowship, during which time he became fluent in German.
In 1960, he took up a lectureship in Newcastle upon Tyne and
in 1965 joined the newly founded Mathematics Institute at Warwick University as a Reader, as part of the original algebra group under the
leadership of Professor Sandy Green. Margaret Green, Sandy Green's widow, recalls that she and Sandy:
\begin{displayquote}{\small ``moved to Warwick in 1965, and Sandy invited
Stewart Stonehewer and Trevor Hawkes and Roger to join him.
As you will know, the Maths Institute building was first of all the library. Christopher Zeeman had a great interest in architecture and design and
created a department like no other; with its huge social area and blackboards! He had obtained large funds with which to bring mathematicians and their
families from around the world. These were exciting and productive times mathematically, and I think it was a very happy time for Roger.''}
\end{displayquote}

Roger was promoted to professor in 1970 and remained at Warwick until his retirement in  2001.
Following his retirement, he continued to live in Leamington Spa until 2019, when he moved to a care home in Wirral, where he died on 21 February 2022.

Several colleagues related the following anecdote, a story which Roger himself often told and which became part of the Warwick local folklore; Charu Hajarnavis
tells the story here:
\begin{displayquote}{\small ``One fine morning Roger came to me in the common room and said, ``You know what I saw this morning or rather didn't see. My garage door was open and 
my car was missing.'' There was no sign of a break in. My first thought was, ``Roger, you parked the car somewhere and returned by an alternative method.''
 I was new and regrettably did not say it. He reported it to the police later that day. It took the police more than a week to find the car parked right
 in front of the local police station within walking distance of Roger's flat. The police were not amused. But they did see the humour after
 they realised that he was a professor of mathematics.''}
\end{displayquote}\bigbreak

\section{His career at Warwick}


\centerline{\it {Roger Carter, the colleague}}\bigbreak

\noindent  Shortly after Roger joined the Mathematics Institute at Warwick, 
 Sandy Green suffered a major stroke and for some time Roger took over informally as the ``head of the algebra group''. Roger served as head of
the mathematics department for four years between 1976 and 1979 and again for one year in 1995.
Roger also co-organised three Warwick thematic symposia, {\it{Finite Groups}}, in 1966/67 with Sandy Green,
{\it{Group Theory}} in 1972/73, and {\it{Groups, Rings and Representations}} in 1990/91. 
In addition, he co-organised two Durham symposia, the first, {\it{Finite Simple Groups}}, in 
1978 with Michael Collins, and the second, {\it{Representation Theory of Algebraic Groups and Related Finite Groups}}, in
1987 with Steve Donkin. In 1997, Roger co-organised, with Michel Brou\'e and Jan Saxl, one of the first  thematic semester-long programmes at the
Isaac Newton Institute for Mathematical Sciences in Cambridge, 
entitled {\it{Representation Theory of Algebraic Groups and Related Finite Groups}}. All of these events attracted a very international audience and contributed to strengthening algebra and group theory in
the UK and at
Warwick in particular.  Michael Collins recalls, for example, that Roger was responsible for inviting Jacques Tits to participate in the 1978 Durham Symposium.
\bigbreak

\noindent Ian Stewart took up a position at Warwick after obtaining his PhD there in 1969. He writes about Roger as a colleague:
\begin{displayquote} {\small 
``Roger was a perfectionist, very much in the school of Carl Friedrich Gauss 
whose motto was \emph{pauca sed matura} (‘few but ripe’). Early in my career
he quietly took me on one side and suggested I was writing too many papers 
too quickly. (This was before the days when administrators assiduously 
counted how many papers you produced every year.)

On one occasion, as a very junior lecturer, I had been put in charge of organising Open Days for prospective students and their parents, and was having trouble finding volunteers to talk to them. I mentioned this to Roger, who promptly got up from his desk and showed me how to do it, by walking briskly round the corridors of the Department banging on office doors and telling people they had just volunteered.''}
\end{displayquote}

\noindent Miles Reid, who joined the Warwick Mathematics Institute in 1978, writes:
\begin{displayquote}{\small ``Roger Carter was chairman at Warwick when I
was hired, and he played a significant role in my
career, giving me very useful advice on the promotion
procedures and on exam administration.''}
\end{displayquote}

\noindent Derek Holt, who also joined the Institute in 1978, speaks of Roger's generous spirit, something echoed by other colleagues and his PhD students:
\begin{displayquote} {\small 
``Roger provided strong and effective leadership of the algebra group at
Warwick for many years, offering help and advice to members of the group
when required. This often involved inviting the person involved for a
meal at a restaurant. I was promoted to professor in 2001, which I think
was shortly after Roger's retirement, but he invited me for a lavish meal
and offered really useful advice on all aspects of my new role.

As Chair he was effective
and efficient and was generally good at arguing
the case for Mathematics at university committees.

In personality he was generally calm and stated his opinions clearly.
I would not exactly describe him as unemotional, because he could  speak
very forcefully on academic topics on which he held strong opinions.''}
\end{displayquote}

\bigbreak


\centerline{\it {Roger Carter, the teacher}}\bigbreak
\begin{displayquote}{\small ``He was a popular lecturer, renowned for giving very clear lectures and well-organised notes, both at undergraduate and postgraduate level.''\quad Ian Stewart} 
\end{displayquote}

\bigbreak
\noindent Anyone who had the good fortune to attend a course, or a colloquium or conference talk, or even a contribution to a working group, prepared and
delivered by Roger, will confirm that Roger's  reputation as an outstanding expositor of mathematics is well-deserved. In addition to a speaking style
which was clear and optimally-paced, Roger's impeccably upright and even hand-writing, whether on the board or on paper, rendered 
discovering mathematics through his presentations an exceptionally enriching experience. Perhaps the following  lines from Derek Holt explain some of his success:
\begin{displayquote} {\small 
``Early in my career I was obliged to go to a number of
generally boring and unhelpful training sessions on things like effective
lecturing, supervising and tutoring, but the single most useful advice on
this topic was a brief remark by Roger: \emph{It is essential that you like the
students attending your lecture. If you like them they will like you, and
if you don't like them then they will not like you}.''}
\end{displayquote}

\noindent Roger's research students speak of his accessibility and supportive supervision:
\begin{displayquote}{\small ``Throughout the process, he met [me] each week, listened
        patiently and carefully, and was encouraging when I finally made sufficient progress. 
I remember one week when he was head of the Institute and the supervision time was in the middle of some crisis where there were people in and
        out of his office in a stressed state. I asked him if he wanted to postpone the meeting but he didn't and switched his mind to it for the duration of
        the meeting. That's the kind of thing you take away from a role model for later use.''\qquad James Franklin}\end{displayquote}\

    \noindent The following story, related to us by Jian-Yi Shi, one of Roger's PhD students, perfectly illustrates Roger's kind and generous nature:
\begin{displayquote}{\small ``Under the introduction of Carter, I met Professor A.O.~Morris in the annual British
Mathematical Colloquium. 
Morris invited me to visit in Aberystwyth University, Wales, in July of
1982. However, the railroad workers went on strike during that period. Thus I could
not go there via railway. In order to support my visit, Carter drove me by his car from
Warwick University to a coach station in the south (it took more than two hours on the
road) where I could take a coach directly to Aberystwyth University. One week later
Carter picked me back to Warwick University from that coach station.''}\end{displayquote}\bigbreak

\bigbreak

\noindent List of Roger Carter's PhD students\bigbreak\begin{enumerate}[]
  \item Mary Jane Prentice (1968) {\it{Normalizers and Covering Subgroups of Finite Soluble Groups}}
\item Gordon Elkington (1970) {\it{Centralisers of Unipotent Elements in Semisimple Algebraic Groups}}
\item Stephen Mayer (1971) {\it{On The Irreducible Characters Of The Weyl Groups}}
\item Philip Gager (1973) {\it{Maximal Tori in Finite Groups of Lie Type}}
\item Stephen Reid (1974) {\it{The Complex Characters of the Adjoint Chevalley Group $C_2(q)$, $q$ Odd}}
\item Pamela Bromwich (1975) {\it{Variations on a Theme of Solomon}}
\item Demetris Deriziotis (1977) {\it{The Brauer Complex and Its Applications to the Chevalley Groups}}
\item James Franklin (1981) {\it{Homomorphisms between Verma Modules and Weyl Modules in Characteristic $p$}}
\item Jian-Yi Shi (1984) {\it{The Decomposition into Cells of the Affine Weyl Groups of Type A}}
\item Alexandros P. Fakiolas (1987) {\it{ Hecke Algebras and the Lusztig Isomorphism}}
\item Mauro Costantini (1989) {\it{On the Lattice Automorphisms of Certain Algebraic Groups}}
\item Maria Iano (1990) {\it{ Polynomial Representations of Symplectic Groups}}
\item Carlos Andr\'e (1992) {\it{Irreducible Characters of the Unitriangular Group and Coadjoint Orbits}}
\item Jacqui Ramagge (1992) {\it{On Some Twisted Kac-Moody Groups}}
\item John Cockerton (1995) {\it{ Canonical Bases and Piecewise-Linear Combinatorics}}
\item Richard Green (1995) {\it{$q$-Schur Algebras and Quantized Enveloping Algebras}}
\item Bethany Marsh (1995) {\it{Canonical Bases and Related Bases in Modules for Quantized Enveloping Algebras}}
  \item Anthony Scott (1996) {\it{Classes of Maximal-Length Reduced Words in Coxeter Groups}}
  \end{enumerate}

\bigbreak

\section{Soluble groups}\label{sec-solvablegrps}

Roger Carter started his research career under the supervision of Derek Taunt, who had himself written a
thesis under the
direction of Philip Hall.
It appears that ``soluble groups'' was the dominant topic at that time among the group theorists in Cambridge.
In Roger's first published article \cite{MR0114859}, coming out of his thesis work, he explains that his work is based upon the theory of
soluble groups developed by Hall.

In order to describe this work, let us introduce some definitions.
Hall had shown that a finite group $G$ is soluble if and only if, for each prime
$p$ dividing its order, $G$ contains a Sylow $p$-complement. One chooses a complete set of
Sylow $p$-complements of a soluble group $G$, one for each prime $p$, and the
set of all intersections of these complements forms, together with $G$, a
lattice of subgroups, called a \emph{Sylow system of $G$}. 
The intersection of the normalisers of the subgroups of a Sylow system 
is called the \emph{normaliser of the system}. The normalisers of the various
Sylow systems of $G$ are called \emph{system normalisers of $G$}. They are all nilpotent
and form a characteristic class of conjugate subgroups of $G$.

In \cite{MR0114859}, Roger defines the notion of an \emph{abnormal} subgroup (a subgroup $H$ of a group $G$ is
abnormal if $g\in\langle H, g^{-1}Hg\rangle$ for all $g\in G$). He then shows that every abnormal subgroup of a soluble
group contains a system normaliser and further that if a system normaliser is equal to its own normaliser,
then it must be abnormal. This then motivates the study of the system normalisers which are equal to their own normalisers.
This class of groups contains, for example, all soluble groups of Fitting height 2.
A second article \cite{MR0133368}, also coming from his thesis work, establishes a connection between the theory of system normalisers and results of
Schenkman \cite{MR0068538} and Higman \cite{MR0079001} on the existence of complements to certain terms in
the lower central series of a soluble group $G$.

In \cite{MR0123603}, Roger established the existence and conjugacy of what would later become known as \emph{Carter subgroups}:
every finite soluble group possesses nilpotent subgroups which are their own normalisers, and all such subgroups are conjugate.
According to Huppert \cite[p. 737]{MR0224703}, the discovery of the Carter subgroups was inspired by the
theory of Lie algebras~---~and Roger certainly was familiar with that theory. For example, in the
introduction of his survey article \cite{MR0174655}, he explicitly acknowledges his indebtedness to an
inspiring course of lectures on Lie algebras given by Hall in Cambridge in 1959.
Roger wrote \cite{MR0123603} during his year in T\"{u}bingen and he thanks Professors Wielandt and Hall for their help
``in the formulation of the main theorem in its present general form and in the construction of the proof''.
In addition, he shows that these subgroups are abnormal and further notes that the nilpotent
self-normalising subgroups satisfy a ``covering and avoidance'' theorem similar to that proven by Hall for
system normalisers (see \cite{MR1575659} for the relevant definitions). In \cite{MR0142652}, he goes on to consider
nilpotent self-normalising subgroups $H$ in a general (not necessarily soluble) finite group $G$ and proves that under two
additional conditions on $H$, namely that $H$ is a Hall subgroup of $G$ and that its Sylow subgroups are regular, $H$
has a normal
complement in $G$. Moreover, he produces examples showing that each of the four given conditions is necessary.

In 1962, Gasch\"{u}tz generalised the notions introduced by Roger, developing the theory of formations (see \cite{MR0179257}).
He introduces the article by saying:
\begin{displayquote}{\small ''The results of the present work have their origin in a
    closer examination of the method Carter used to construct self-normalising nilpotent subgroups in finite
    soluble groups.''}\end{displayquote}
A class $\mathfrak{F}$  of finite  groups is called a {\it formation} if (1) homomorphic images
of groups in  $\mathfrak{F}$ are also in  $\mathfrak{F}$, and (2) if $N_1$ and $N_2$ are normal subgroups of a group $G$ such that $G/N_1 \in  \mathfrak{F}$ and $G/N_2 \in  \mathfrak{F}$, then $G/(N_1 \cap N_2) \in  \mathfrak{F}$.
Gasch\"{u}tz further introduces the concept of a {\it saturated} formation to be a formation  $\mathfrak{F}$ required  to contain certain extensions of groups
in $\mathfrak{F}$ (see \cite[Kapitel VI, \S7]{MR0224703}).  The main result of \cite{MR0179257} states that, for a saturated formation  $\mathfrak{F}$, each finite
soluble group $G$ possesses an $\mathfrak{F}$-{\it covering subgroup} $H$ (see \cite{MR0224703});
furthermore,  all such are conjugate in $G$. Since the class of nilpotent groups $\mathfrak{N}$ is a saturated formation and self-normalising nilpotent subgroups of a
soluble group $G$ are  $\mathfrak{N}$-covering subgroups, the existence and conjugacy of Carter subgroups becomes a corollary of Gasch\"{u}tz's result.

Roger relates his ``lucky'' meeting with Gasch\"{u}tz at that time, \cite{zbMATH05504621}:

\begin{displayquote}{\small ``In 1960 I took up a lectureship in Newcastle and in 1962 I paid a visit to Stockholm to attend the International Congress of Mathematicians.
    During this meeting I learned from Joachim Neub\"{u}ser that Gasch\"{u}tz had obtained new results on finite soluble groups containing my results as a special case.
    As the train back from Stockholm to London stopped quite close to Kiel, I decided to take the opportunity to get out and visit Professor Gasch\"{u}tz.
    When I arrived at Kiel I found that he was at a family birthday party, but he very kindly arranged to have a discussion with me, and told me about his new theory of formations and covering subgroups.
    This conversation inspired me to start working on formations and the so-called ${\mathfrak {F}}$-normalisers in a finite soluble group, which generalise the system normalisers which I had studied in my PhD thesis.
    I subsequently found that Trevor Hawkes had obtained similar results independently, and so we published them in a joint paper, which eventually appeared in 1967.'' (See \cite{MR0206089}.) }\end{displayquote}

It seems that Gasch\"{u}tz's work was the driving force for the further study of finite soluble groups at that time.  Roger
went on to collaborate on one further article \cite{MR0228581} on soluble groups. Trevor Hawkes relates:

\begin{displayquote}{\small ``In 1966 I got an invitation from Reinhold Baer to spend a couple of weeks in Frankfurt and met Bernd Fischer there. Roger was also in contact with
    Fischer, who was interested in finite soluble groups at that time.The following year I left Cambridge to join the new Mathematics Department at the recently-founded University
    of Warwick. Roger Carter and Sandy Green were already there, and Brian Hartley subsequently joined us, initially commuting from Cambridge on a motorbike. While Fischer was on an extended visit to Warwick, he, Roger and I  wrote a joint paper, called ``Extreme Classes of Finite Soluble Groups'' \cite{MR0228581}. Those were
    exciting times when the Mathematics Institute was in an isolated house on Gibbet Hill, more than a mile from the central Warwick Campus, and my office was a former bathroom.
    Sandy Green, who had worked at Bletchley Park along with many other mathematicians during the war (and had met his wife Margaret there), was head of our Algebra Group.  Roger took over temporarily during Sandy's recovery from a stroke that kept him out of action for a year or more.''}
\end{displayquote}

While Hawkes continued to work on soluble groups, producing the massive volume \cite{MR1169099}, joint with Klaus Doerk,
by the late 1960's Roger had moved on to the groups of Lie type. The first signs of that
transition already became visible as early as 1963, when Roger gave a lecture on ``Simple groups
and simple Lie algebras'' at the British Mathematical Colloquium at the R.M.C. Shrivenham (UK). Shortly
afterwards appeared Roger's joint article with Paul Fong \cite{MR0166271}, his third most-often cited article: 
\begin{displayquote}{\small ``Roger Carter and I first met in Stockholm at the 1962 Congress.
    By good luck we were both interested then in the Sylow 2-subgroups of classical groups in odd characteristic.''
  \qquad P.~Fong}\end{displayquote} In their article,
the authors
determine the Sylow $2$-subgroups and their normalisers for all families of finite classical groups defined over a field of
odd characteristic. The stage had been set for their work by the paper of A.~Weir \cite{MR0072143} which handled the Sylow $l$-subgroups of the
finite classical groups defined over a field of characteristic $p\ne l$, for odd primes~$l$. A look at the MathSciNet citation list shows that the results of the article
are still regularly used in
a wide range of research topics, including character theory,
modular representation theory,
fusion systems, and subgroup structure. The theory was extended when various researchers turned to the Sylow subgroups of the
finite simple groups of Lie type. For a discussion of the further developments and the most ``modern'' treatment, one can refer to \S25.2 of \cite{MR2850737}.
Note however that the Carter--Fong paper remains the definitive work for the question they considered and solved.

\section{Conjugacy classes} \label{sec-conjclasses}

Let $G$ be a simple algebraic group over an algebraically closed field
$k$ of characteristic~$p$. (Initially, $p$ may be zero or a prime; later
on, we will assume that $p>0$.) We refer to \cite{MR2850737} (or the
first chapters of \cite{MR1266626}) for general background on algebraic 
groups. Let $B$ be a Borel subgroup of $G$ and $T$ be a maximal torus
contained in $B$. Then $W:=N_G(T)/T$ is the Weyl group of $G$ and we have 
the Bruhat decomposition
\[ G=\coprod_{w \in W} B\dot{w}B \qquad \mbox{(where
$\dot{w}$ is a representative of $w$ in $N_G(T)$)}.\]
This is one of the main structural properties of $G$. It is a key tool
which allows us to reduce many questions about $G$ to questions about
the Weyl group $W$. On a basic level, this already happens when we consider 
the conjugacy classes of $G$ and of $W$. On a more elaborate level, this 
continues when we consider representations of $G$ and of related finite
groups. Let us now focus on conjugacy classes. Roger Carter made 
substantial, and highly influential contributions to this circle of ideas; 
the references are (in chronological order): \cite{MR0269749}, 
\cite{MR0318337}, \cite{MR0289521}, 
\cite{MR0417306}, \cite{MR0417307}, 
\cite{MR0512022}, \cite{MR0602121}, \cite{MR1973577}. 

\begin{table}[htbp] 
\caption{Dynkin diagrams of finite type} \label{Mdynkintbl} 
\begin{center} {\small 
\begin{picture}(290,155)
\put( 00, 25){$E_7$}
\put( 20, 25){\circle*{5}}
\put( 18, 30){$1$}
\put( 20, 25){\line(1,0){20}}
\put( 40, 25){\circle*{5}}
\put( 38, 30){$3$}
\put( 40, 25){\line(1,0){20}}
\put( 60, 25){\circle*{5}}
\put( 58, 30){$4$}
\put( 60, 25){\line(0,-1){20}}
\put( 60, 05){\circle*{5}}
\put( 65, 03){$2$}
\put( 60, 25){\line(1,0){20}}
\put( 80, 25){\circle*{5}}
\put( 78, 30){$5$}
\put( 80, 25){\line(1,0){20}}
\put(100, 25){\circle*{5}}
\put( 98, 30){$6$}
\put(100, 25){\line(1,0){20}}
\put(120, 25){\circle*{5}}
\put(118, 30){$7$}

\put(150, 25){$E_8$}
\put(170, 25){\circle*{5}}
\put(168, 30){$1$}
\put(170, 25){\line(1,0){20}}
\put(190, 25){\circle*{5}}
\put(188, 30){$3$}
\put(190, 25){\line(1,0){20}}
\put(210, 25){\circle*{5}}
\put(208, 30){$4$}
\put(210, 25){\line(0,-1){20}}
\put(210, 05){\circle*{5}}
\put(215, 03){$2$}
\put(210, 25){\line(1,0){20}}
\put(230, 25){\circle*{5}}
\put(228, 30){$5$}
\put(230, 25){\line(1,0){20}}
\put(250, 25){\circle*{5}}
\put(248, 30){$6$}
\put(250, 25){\line(1,0){20}}
\put(270, 25){\circle*{5}}
\put(268, 30){$7$}
\put(270, 25){\line(1,0){20}}
\put(290, 25){\circle*{5}}
\put(288, 30){$8$}

\put( 00, 59){$G_2$}
\put( 20, 60){\circle*{6}}
\put( 18, 66){$1$}
\put( 20, 58){\line(1,0){20}}
\put( 20, 60){\line(1,0){20}}
\put( 20, 62){\line(1,0){20}}
\put( 26, 57.5){$>$}
\put( 40, 60){\circle*{6}}
\put( 38, 66){$2$}

\put( 80, 60){$F_4$}
\put(100, 60){\circle*{5}}
\put( 98, 65){$1$}
\put(100, 60){\line(1,0){20}}
\put(120, 60){\circle*{5}}
\put(118, 65){$2$}
\put(120, 58){\line(1,0){20}}
\put(120, 62){\line(1,0){20}}
\put(126, 57.5){$>$}
\put(140, 60){\circle*{5}}
\put(138, 65){$3$}
\put(140, 60){\line(1,0){20}}
\put(160, 60){\circle*{5}}
\put(158, 65){$4$}

\put(190, 75){$E_6$}
\put(210, 75){\circle*{5}}
\put(208, 80){$1$}
\put(210, 75){\line(1,0){20}}
\put(230, 75){\circle*{5}}
\put(228, 80){$3$}
\put(230, 75){\line(1,0){20}}
\put(250, 75){\circle*{5}}
\put(248, 80){$4$}
\put(250, 75){\line(0,-1){20}}
\put(250, 55){\circle*{5}}
\put(255, 53){$2$}
\put(250, 75){\line(1,0){20}}
\put(270, 75){\circle*{5}}
\put(268, 80){$5$}
\put(270, 75){\line(1,0){20}}
\put(290, 75){\circle*{5}}
\put(288, 80){$6$}

\put( 00,110){$D_n$}
\put( 00,100){$\scriptstyle{n \geq 3}$}
\put( 20,130){\circle*{5}}
\put( 25,130){$1$}
\put( 20,130){\line(1,-1){21}}
\put( 20, 90){\circle*{5}}
\put( 26, 85){$2$}
\put( 20, 90){\line(1,1){21}}
\put( 40,110){\circle*{5}}
\put( 38,115){$3$}
\put( 40,110){\line(1,0){30}}
\put( 60,110){\circle*{5}}
\put( 58,115){$4$}
\put( 80,110){\circle*{1}}
\put( 90,110){\circle*{1}}
\put(100,110){\circle*{1}}
\put(110,110){\line(1,0){10}}
\put(120,110){\circle*{5}}
\put(117,115){$n$}

\put(170,110){$C_n$}
\put(170,100){$\scriptstyle{n \geq 2}$}
\put(190,110){\circle*{5}}
\put(188,115){$1$}
\put(190,108){\line(1,0){20}}
\put(190,112){\line(1,0){20}}
\put(196,107.5){$>$}
\put(210,110){\circle*{5}}
\put(208,115){$2$}
\put(210,110){\line(1,0){30}}
\put(230,110){\circle*{5}}
\put(228,115){$3$}
\put(250,110){\circle*{1}}
\put(260,110){\circle*{1}}
\put(270,110){\circle*{1}}
\put(280,110){\line(1,0){10}}
\put(290,110){\circle*{5}}
\put(287,115){$n$}

\put( 00,150){$A_n$}
\put( 00,140){$\scriptstyle{n \geq 1}$}
\put( 20,150){\circle*{5}}
\put( 18,155){$1$}
\put( 20,150){\line(1,0){20}}
\put( 40,150){\circle*{5}}
\put( 38,155){$2$}
\put( 40,150){\line(1,0){30}}
\put( 60,150){\circle*{5}}
\put( 58,155){$3$}
\put( 80,150){\circle*{1}}
\put( 90,150){\circle*{1}}
\put(100,150){\circle*{1}}
\put(110,150){\line(1,0){10}}
\put(120,150){\circle*{5}}
\put(117,155){$n$}

\put(170,150){$B_n$}
\put(170,140){$\scriptstyle{n \geq 2}$}
\put(190,150){\circle*{5}}
\put(188,155){$1$}
\put(190,148){\line(1,0){20}}
\put(190,152){\line(1,0){20}}
\put(196,147.5){$<$}
\put(210,150){\circle*{5}}
\put(208,155){$2$}
\put(210,150){\line(1,0){30}}
\put(230,150){\circle*{5}}
\put(228,155){$3$}
\put(250,150){\circle*{1}}
\put(260,150){\circle*{1}}
\put(270,150){\circle*{1}}
\put(280,150){\line(1,0){10}}
\put(290,150){\circle*{5}}
\put(287,155){$n$}
\end{picture}}
\end{center}
\begin{center} 
{\footnotesize (The numbers attached to the vertices define 
a standard labelling of the graph.)}
\end{center}
\end{table}

By Chevalley's famous S\'eminaire \cite{MR2124841}, the simple
algebraic groups $G$ as above are classified by the Dynkin diagrams in 
Table~\ref{Mdynkintbl} (plus some extra information about the root and 
weight lattices which we can ignore in the discussion below). Thus, it
is often possible to obtain results about $G$ based on ``case--by--case'' proofs.
As far as the Weyl group~$W$ is concerned~---~and 
prior to Carter \cite{MR0269749}, \cite{MR0318337}~---~the conjugacy 
classes had been determined individually in all cases. For example, $W$ of 
type $A_{n-1}$ ($n\geq 2$) is isomorphic to the symmetric group $\fS_n$, 
and it is well-known that the conjugacy classes are parametrised 
by the partitions of~$n$. For the other infinite series $B_n$, $C_n$,
$D_n$ similar descriptions in terms of combinatorial objects are available. 
For $G$ of exceptional type, explicit lists of class representatives, 
centraliser orders etc.\ are known in the form of tables; e.g., for~$W$ of 
type~$E_8$, there are precisely $112$ conjugacy classes. Nowadays, all this 
is readily available in electronic databases; e.g., within the {\sf CHEVIE}
system. (See J.~Michel's enlarged, and more recent version of the 
original system at \url{https://github.com/jmichel7/Chevie.jl}.) 
Nevertheless, it is certainly desirable to exhibit unifying principles and 
general patterns behind such ``case--by--case'' classifications. Carter's 
work \cite{MR0269749}, \cite{MR0318337} achieves exactly that: the 
conjugacy classes of $W$ are classified in terms of~what he called 
``admissible graphs'', which can be seen as generalisations of the 
partitions of $n$ which are used for type $A_{n-1}$. In particular, it
is shown that every element $w\in W$ can be written as $w=uv$ where $u,v
\in W$ are involutions; so $w$ is conjugate to its inverse $w^{-1}$. 
These results are still the best reference in this area. 

Now consider the conjugacy classes of $G$. By the Jordan decomposition
of elements, any $g\in G$ can be written uniquely as $g=su=us$ where
$s\in G$ is semisimple and $u\in G$ is unipotent. Consequently, the 
classification of conjugacy classes of $G$ is reduced to two separate 
problems: 1)~classify the conjugacy classes of semisimple elements and 
2)~classify the conjugacy classes of unipotent elements in the centraliser
$C_G(s)$, for every semisimple $s\in G$. Step~1) is by now well-understood; 
see, e.g., the work of Carter's student D.~Deriziotis \cite{MR0742140}. 
(G.~Lusztig pointed out to us that the classification of semisimple classes
actually goes back to V. G. Kac \cite{MR0251091}, and even further to the 
work of de Siebenthal in the 1930s.) Carter himself made substantial 
contributions to the description of the centralisers of semisimple elements;
see \cite{MR0512022}, \cite{MR0602121}. Such information is not merely of 
interest in its own right but it is also relevant to the character theory of 
finite Chevalley groups (via Lusztig's Jordan decomposition of characters). 

Let us now look at the case $s=1$, that is, the problem of classifying the 
conjugacy classes of unipotent elements of~$G$ itself. The story is 
more involved since~$p$, the characteristic of the field over which 
$G$ is defined, plays a role here. One important general result, at least,
is the fact that there are only finitely many unipotent classes; see
Lusztig \cite{MR0419635} and further references there. A classification 
had been known for $p=0$ since the 1960's by a general scheme (``weighted
Dynkin diagrams''). A detailed exposition of the methods required to prove 
this can be found in Chapter~5 of Carter's book \cite{MR1266626}; actually, 
this also works when $p$ is a prime that is larger than an explicitly 
given lower bound. In \cite{MR0417306}, \cite{MR0417307}, Carter 
and his masters student P.~Bala found a new, and particularly simple 
description of the unipotent classes of~$G$, by using the concept of
``distinguished'' parabolic subgroups; again, this works for $p=0$ 
or~$p$ sufficiently large. Eventually, it was proved that 
\begin{center}
\textit{the Bala--Carter classification always works unless~$p$ is a 
``bad'' prime for~$G$};
\end{center}
see \cite[\S 5.11]{MR1266626} and further references there. (A prime $p$ 
is called ``bad'' if $p=2$ for $G$ of types $B_n$, $C_n$, $D_n$; $p=2$
or~$3$ for types $G_2$, $F_4$, $E_6$, $E_7$; $p=2$,~$3$~or~$5$ for type 
$E_8$.) These results have been widely used, and they provide a standard 
nomenclature for unipotent classes; see, for example, Liebeck--Seitz 
\cite{MR2883501}. A more modern approach to the Bala--Carter 
classification~---~which uniformly works for any $p$ that is not bad~---~is 
due to Premet \cite{MR1976699}. But things are definitely different in the 
bad prime case; see \cite{MR2883501} for further references and explicit 
classifications in all cases, including $p$ bad.

However, and quite remarkably, the story does not end here for Bala--Carter.
Indeed, Spaltenstein \cite[\S 5.2]{MR0672610} shows that one can always 
naturally embed the set of unipotent classes of $G_0$ (where $G_0$ is a 
group of the same type as $G$ but defined over the field of complex 
numbers) into the set of unipotent classes of~$G$. Thus, one may speak of 
those unipotent classes of $G$ that ``come from characteristic~$0$''.
Furthermore, there is a series of articles by Lusztig, starting with 
\cite{MR2183120}, which presents a picture of unipotent classes of $G$ 
which should apply to any~$p$ and is as 
close as possible to the known picture in characteristic~$0$. The main idea 
in Lusztig's papers is to define (by a general scheme) a partition of the 
unipotent variety of~$G$ into finitely many locally closed and $G$-stable
``pieces'', which are in bijection with the unipotent classes in 
characteristic~$0$. If $p$ is not a bad prime, then each piece
consists of a single unipotent class. Otherwise, each piece is the 
union of a unique unipotent class that comes from characteristic~$0$, and 
possibly a (small) number of ``new'' unipotent classes. Thus, the 
Bala--Carter classification can still be used to parametrise the 
pieces of the unipotent variety. The picture is illustrated by 
Table~\ref{MuniG2} for $G$ of type~$G_2$.


\begin{table}[htbp] 
\caption{Unipotent classes of $G$ of type $G_2$} \label{MuniG2} 
\begin{center} {\small $\begin{array}{c@{\hspace{2cm}}c@{\hspace{2cm}}c} 
\hline \text{Class $\cO$} & \text{condition on $p$} & \dim \cO \\\hline 
\varnothing & \text{any $p$} &0 
\\ A_1 & \text{any $p$} &6 \\
(\tilde{A}_1)_{3} & \text{only $p=3$} & 6 \\
\tilde{A}_1 & \text{any $p$} &8 \\
G_2(a_1) & \text{any $p$} &10 \\ 
G_2 & \text{any $p$} &12 \\ \hline 
\multicolumn{3}{c}{\text{(See Lawther \cite[Table~B]{MR1351124} or 
Spaltenstein \cite[\S 5]{MR0803340})}}
\end{array}$}
\end{center}
\end{table}


In that table we see that, from characteristic~$0$, we always obtain $5$ 
unipotent classes denoted by $\varnothing$, $A_1$, $\tilde{A}_1$, $G_2(a_1)$,
$G_2$, exactly as in the table on p.~8 of Bala--Carter \cite{MR0417307}; 
when $p=3$, there is a ``new'' unipotent class denoted by $(\tilde{A}_1)_3$. 
The classes $\tilde{A}_1$ and $(\tilde{A}_1)_3$ form a piece in the 
sense of Lusztig; all other unipotent classes form a piece by themselves. 


Finally we mention a more recent development concerning a relationship 
between conjugacy classes in $W$ and unipotent classes of $G$. On p.~95 
in \cite{MR0865891}, Carter reports on the results of the PhD thesis 
of A.~J.~Starkey under the direction of his colleague Sandy Green. These 
results show that elements of minimal length in the conjugacy classes
of $W=\fS_n$ have some special properties with respect to representations of 
the Iwahori--Hecke algebra associated with $\fS_n$. Analogous results for 
$W$ of any type have subsequently been obtained in \cite{MR1250466}. This is 
one ingredient in Lusztig's construction \cite{MR2833465} of a canonical 
surjective map 
\[ \Phi\colon \{\text{conjugacy classes of $W$}\} \rightarrow 
\{\text{unipotent classes of $G$}\}.\]
This essentially relies on studying the intersections $\cO \cap B\dot{w}B$ 
where $\cO$ is a unipotent class of $G$ and $w\in W$ is an element that has 
minimal length in its conjugacy class in $W$. With this result of 
Lusztig, we return full circle to the initial remarks in this section 
about the importance of the Bruhat decomposition. And the idea that a 
relationship between conjugacy classes of~$W$ and unipotent classes of 
$G$ might exist already appeared in the final section of Carter's work 
\cite{MR0318337}; see also \cite{MR0289521}. In a somewhat different
direction, the problem of finding the cardinalities of intersections
$C\cap B\dot{w}B$, where $C$ is a conjugacy class in $\mbox{GL}_n(q)$, 
came up in Carter's paper \cite{MR1152995}.

\section{Modular representations} \label{sec-modular}

Roger Carter made several important contributions (partly in joint work) 
to the modular representation theory of the symmetric group $\fS_n$, of 
the general linear group $\mbox{GL}_n$ and, more generally, of a simple 
algebraic group $G$ (as in the previous section) and related finite groups. 
By ``modular'' we mean representations over fields of positive 
characteristic $p>0$; this is often much harder than working over fields 
of characteristic~$0$ (where, in the above situations, representations are 
semisimple). In each case, one of the main problems is to determine the 
simple (or irreducible) modules.

Let us begin by looking at $\fS_n$. By the classical work of Frobenius 
we know that the simple $\fS_n$-modules over $\CC$ are parametrised by 
the partitions $\lambda\vdash n$. Now let $k$ be a field of characteristic
$p>0$. The fundamental work of G. D. James in the 1970s (see, e.g., the 
exposition by A. Mathas \cite{Mathas}) leads to a parametrisation of the 
simple $\fS_n$-modules over~$k$, in two steps. Firstly, there are the 
``Specht modules'' $\{S^\lambda \mid \lambda \vdash n$\} which are defined 
and free over $\ZZ$, with standard bases given by certain combinatorial 
objects (Young tableaux); they remain irreducible when we extend scalars 
from $\ZZ$ to $\CC$, and all simple $\fS_n$-modules over $\CC$ arise in this 
way. On the other hand, if we extend scalars from $\ZZ$ to $k$, then 
$S^\lambda$ may become reducible. Now James showed that if 
$\lambda\vdash n$ is $p$-regular, then the $p$-modular reduction of 
$S^\lambda$ has a unique simple quotient, called $D^\lambda$; furthermore,
in this way, one obtains all simple $\fS_n$-modules over $k$. Note, however,
that finding an explicit formula for $\dim D^\lambda$ is still an open 
problem! (See \cite[\S 5.9]{Mathas} for further comments on this.)

Roger was interested in the following natural question: Which Specht modules
$S^\lambda$, for $p$-regular $\lambda\vdash n$, actually remain irreducible 
modulo~$p$? He conjectured an answer to this question, in terms of
a precise purely combinatorial condition on the partition~$\lambda$. 
He never published his conjecture, but the subject was taken up~---~and 
completely solved~---~in the articles by James \cite{MR0463281} and 
James--Murphy \cite{MR0541676}; see \cite[\S 5.4]{Mathas} for further details.

Carter's article \cite{MR0933425} on raising and lowering operators for 
the special linear Lie algebra~$\mathfrak{sl}_n$ was not on modular 
representation theory as such. But it inspired the definition of 
``modular lowering operators'' in the fundamental work of A. Kleshchev 
\cite{MR1319521} on modular branching rules for~$\fS_n$.

Let us now turn again to a simple algebraic group $G$ over an algebraically 
closed field~$k$ of characteristic~$p$. We will be interested in studying the 
finite-dimensional ``rational'' $G$-modules over~$k$; see, e.g., Humphreys 
\cite{MR2199819} for general background. Again, if $p>0$, then the 
determination of the simple $G$-modules has turned out to be a very 
hard problem. Similarly to what we saw in terms of Specht modules for 
$\fS_n$, a key idea in this context is to approach the simple $G$-modules via 
a process of $p$-modular reduction from certain ``universal'' modules defined
over~$\ZZ$. An explicit construction of such modules over $\ZZ$ for $G$ of
type $A_{n-1}$, and the name of ``Weyl modules'' for their $p$-modular 
reductions, appeared in the highly influential 1974 paper \cite{MR0354887} 
by Carter and Lusztig. For the definition of Weyl modules in general, see
Humphreys \cite[\S 3.1]{MR2199819}. The submodule structure of the Weyl 
modules in the simplest possible case (type $A_1$) is described by 
Carter--Cline \cite{MR0399281}. 

One of the aims of Carter--Lusztig \cite{MR0354887} was to try to extend 
as much as possible the classical Schur--Weyl duality between (polynomial) 
representations of $\mbox{GL}_n$ and representations of the symmetric group 
$\fS_n$ (both over $\CC$), to the case of characteristic $p>0$. Much of the
paper is concerned with the construction of certain non-trivial homomorphisms 
between different Weyl modules associated to certain reflections in the 
affine Weyl group. This confirmed ideas of D.~N.~Verma and played a role 
in the later formulation of a character formula for the simple $G$-modules 
in terms of the affine Weyl group, that is, the fundamental ``Lusztig
Conjecture''; see \cite[\S 3.11]{MR2199819}. 

A substantial part of Carter--Lusztig \cite{MR0354887} has been
integrated into Sandy Green's equally influential volume \cite{MR2349209} 
in the Springer LNM series. See also \cite[p.~10]{MR2349209} for further
historical remarks. In particular, as A.~Kleshchev pointed out to us, there 
are remarkable connections between \cite{MR0354887} and combinatorics and 
invariant theory through the work of D\'esarm\'enien, Kung and Rota in the 
1970s; see M. Clausen \cite{MR0544848}. Since the pioneering work of Carter 
and Lusztig, homomorphisms between Weyl modules and between Specht modules 
have been studied in many papers, including further work \cite{MR0556922} by 
Carter and his student M.~T.~J.~Payne. James Franklin relates the following
amazing story:

\begin{displayquote} 
{\small ``I remember Roger telling me that he (Payne) completed his PhD 
thesis, but posted the only copy of it to Roger and it was lost in the post.
He (Roger) made inquiries with the post office but it was never found.
Presumably the joint paper writes up some results from the thesis.''
\qquad James Franklin}
\end{displayquote}

Finally, let us drop the assumption that $k$ is algebraically closed
of characteristic $p>0$; instead we assume that $k$ is a finite field
of characteristic~$p>0$ and let $G$ be a Chevalley group over $k$ (a
finite group, see Section~\ref{sec-two} below for the general definition). 
Let $U$ be a Sylow $p$-subgroup of $G$. The simple $G$-modules over $k$ 
had been classified by C. W. Curtis and R. Steinberg; a model for them 
in terms of the induced representation $\mbox{Ind}_U^G(1)$ was given by 
C. W. Curtis and F. Richen in the 1960s. (See the expositions in
Curtis--Reiner \cite[\S 72B]{MR0892316}, Humphreys 
\cite[Chap.~7]{MR2199819}, and also further references there.)

Carter and Lusztig \cite{MR0396731} give another model for the 
simple $G$-modules in terms of some explicit endomorphisms of 
$\mbox{Ind}_U^G(1)$ defined using certain elements in the Yokonuma 
algebra (a generalization of the Iwahori-Hecke algebra) with coefficients
in $k$. As Lusztig pointed out to us, it would be interesting 
to use this description to get information on the characters of the 
simple $G$-modules over~$k$.

The two joint papers \cite{MR0354887}, \cite{MR0396731} were written while
Carter and Lusztig were colleagues at Warwick. As Caroline Series remembers:

\begin{displayquote} {\small ``George Lusztig came to Warwick as a Research 
Fellow in 1971 and  left for MIT in 1978, the year I arrived in Warwick.
Although I never met him, it was clear that Roger held him in the highest
esteem and viewed his time in Warwick as having been of great importance. 
Even though I was far from the field, on several occasions years later 
Roger talked to me about him; it was clear this was one of the major 
influences in his mathematical life.''}
\end{displayquote}

%




\section{The two books on groups of Lie type} \label{sec-two}

\begin{displayquote}{\small``Those who have read his books will know that 
he wrote clearly, concisely, and without any hint of pretentiousness. Whether 
in lectures or in writing, Roger's primary aim was always to communicate the 
mathematics in a way that we would understand.''\qquad Jacqui Ramagge}
\end{displayquote}\bigbreak

The two expository books \cite{MR1013112}, \cite{MR1266626} on 
groups of Lie type stand out in Roger's bibliography. The first one 
has more than $900$ citations in MathSciNet; the second one more 
than~$1300$. The two volumes have been highly influential and represent
marvelous pieces of scholarship and service to the general mathematical 
community. There are excellent and highly informative reviews in MathSciNet
written by L.~Solomon (for \cite{MR1013112}) and by D.~B.~Surowski 
(for \cite{MR1266626}); so we will just add some comments.

C. Chevalley published a seminal paper \cite{MR0073602} in 1955, in 
which he introduced the groups that now bear his name. Let $\fg$ be a 
finite-dimensional simple Lie algebra over $\CC$, with root system $R$. 
For each $\alpha\in R$ we have a root element $e_\alpha\in \fg$, which is 
unique up to multiplication by a non-zero scalar. Any choice of a system 
of root elements gives rise to a well-defined basis of $\fg$. Now consider 
the ``divided powers'' $e_\alpha^{[d]}:=\mbox{ad}(e_\alpha)^d/d!\colon 
\fg\rightarrow\fg$ for any $d\geq 0$, where $\mbox{ad}(e_\alpha)(x):=
[e_\alpha,x]$ for $x\in\fg$. Then Chevalley's fundamental insight was that 
the root elements can be chosen such that, with respect to the 
corresponding basis of $\fg$ (now called a ``Chevalley basis''), the 
matrices of the operators $e_\alpha^{[d]}$ have integer coefficients. 
Consequently, if $K$ is an arbitrary field (in particular, a finite field), 
we obtain invertible operators over $K$ by setting $x_\alpha(t):=
\sum_{d\geq 0} t^d e_\alpha^{[d]}$ for any $t\in K$. (Note that $\mbox{ad}
(e_\alpha)$ is nilpotent and so the sum is finite.) The Chevalley group 
over $K$ corresponding to $\fg$ is then defined by $G_K:=\langle 
x_\alpha(t)\mid \alpha\in R,t\in K\rangle$. These groups always have a 
trivial center; in fact, with a few exceptions, $G_K$ turns out to be 
simple. In this way, Chevalley exhibited several families of previously 
unknown simple groups.

In 1965 Roger published a longer survey article \cite{MR0174655} on simple
groups and simple Lie algebras, which was written while he was still at 
Newcastle University. Following the joint work \cite{MR0166271} with Fong 
on Sylow subgroups of classical groups, this was Roger's first paper of
a genuine ``Lie-theoretic'' flavour; it shows that he was already perfectly 
familiar, not only with Chevalley's original work, but also with all the
more recent further developments due to Ree, Steinberg, Suzuki, Tits and 
others. M. Collins pointed out to us that 
\begin{displayquote} {\small ``this survey was a key paper to read as the 
only single source on all the then known finite simple groups, at the time 
that these became of interest to the UK community largely as a consequence 
of the 1966/67 programme in Warwick.''}
\end{displayquote}
Roger's subsequent book \cite{MR1013112} contains an exposition of 
Chevalley's construction, and develops the content of his survey
\cite{MR0174655} in full detail, including order formulae (when $K$ is 
finite), the Steinberg presentation, automorphisms and the twisted
 Chevalley groups. The needed general results about the structure of 
simple Lie algebras over $\CC$ are explained without proofs (but with 
appropriate references); from then on, any other techniques which are used 
in the development are explained (and proved) from first principles. A 
considerable part of the purely group-theoretical aspects is developed 
in the context of groups with a $(B,N)$-pair, as introduced by J.~Tits.

Finally, two further comments. The first one is about a recent development
concerning the above-mentioned ``Chevalley bases'' of $\fg$. Note that 
these are not unique; one is free to individually replace a certain number 
of root elements $e_\alpha$ in such a basis by $\pm e_\alpha$, and one still 
obtains a Chevalley basis. The amount of freedom is described 
in \cite[Prop.~4.2.2]{MR1013112}, in terms of the concept of ``extraspecial
pairs'' of roots. Now, using Lusztig's theory of ``canonical bases'' (to 
be discussed below), it is possible to single out a canonical Chevalley 
basis which is unique up to one global sign; see \cite{MR3652779}, 
\cite{MR3589909}. \label{candj} In particular, this gives rise to a truly
``canonical'' model of $G_K$ for any field~$K$; see the explicit formulae 
in \cite[\S 2.2]{MR3589909}.

The second comment is about the comparison with Steinberg's formidable 
``Lectures on Chevalley groups'' \cite{MR3616493} (which appeared after
Carter's survey \cite{MR0174655} but before the book~\cite{MR1013112}).
Firstly, Steinberg's notes are more densely written, with significantly fewer 
detailed computations.  Furthermore, they also deal with aspects of the 
representation theory of $G_K$, and with Chevalley groups of ``non-adjoint 
type'', which may have a (finite) non-trivial center. All this requires 
further material from the theory of Lie algebras and their representations. 
Moreover, if the ground field $K$ is algebraically closed, then the groups 
$G_K$ are seen to be simple algebraic groups. In Steinberg's more general
setting, one actually obtains constructions for all possible simple 
algebraic groups (as in Section~\ref{sec-conjclasses}). Thus, Carter's
account is somewhat more restrictive but it may be useful in providing a 
gentler introduction, e.g., for graduate students, or for those new to 
the theory of Chevalley groups. As M. Costantini remembers:
\begin{displayquote}
{\small ``I had been studying his book ``Simple groups of Lie type'' for 
my first degree in Italy under the supervision of Giovanni Zacher. It was 
Zacher himself to advise me to take a PhD at Warwick under the supervision 
of Carter, and I could not have had a better advice.''}
\end{displayquote}

Let us now turn to \cite{MR1266626}, the starting point of which is 
the seminal article \cite{MR0393266} by P.~Deligne and G.~Lusztig (1976). 
In that article, the all-important generalised characters $R_{T,\theta}$ 
of a Chevalley group over a finite field are introduced and studied. The 
construction of $R_{T,\theta}$ uses the theory of $\ell$-adic cohomology, 
which was developed by M. Artin, A. Grothendieck et al. in the famous 
``S\'eminaire de G\'eometrie Alg\'ebrique du Bois Marie'' (SGA) that ran 
from about 1960 to 1969 at the IHES near Paris, France. That theory played a 
major role in the solution of the Weil conjectures which, conversely, 
actually had been a driving force for the development of the theory itself.
Carter's book \cite{MR1266626} contains an appendix which
explains the Weil conjectures and gives an introduction to the theory 
of $\ell$-adic cohomology. Given an algebraic variety~$X$ (possibly over 
a field of positive characteristic), that theory produces a family of 
finite-dimensional vector spaces $H_c^i(X)$ ($i\in \ZZ$) over a field of 
characteristic~$0$. Now Deligne--Lusztig define certain varieties $X$ on 
which a Chevalley group $G$ over a finite field acts; then each $H_c^i(X)$ 
naturally becomes a $G$-module and this eventually leads to the definition 
of $R_{T,\theta}$. 

One of the remarkable achievements of Carter's book is
that he succeeds in making this highly elaborate machinery accessible to a 
reader, even with only a modest background in algebraic geometry. More 
specifically, Carter does not assume that the actual 
definition of $H_c^i(X)$ (as explained in the appendix of his book) is 
understood, or even known at all! Instead, in \S 7.1, there is a list of 
formal properties that are satisfied by the cohomology spaces $H_c^i(X)$. 
These properties are then taken as axioms on which are based the subsequent
developments of the representation theory. For example, in the initial theory
of the generalised characters $R_{T,\theta}$ there are two crucial results
(see \cite[\S 6]{MR0393266}): the orthogonality relations (which give 
explicit formulae for the inner product of two generalised characters $R_{T,
\theta}$ and $R_{T',\theta'}$) and the ``exclusion theorem'' (which gives 
a precise condition on when $R_{T,\theta}$ and $R_{T',\theta'}$ do not
have any irreducible constituent in common). The exposition of the proofs 
of these two results on pp.~212--225 in \cite{MR1266626} provides a 
superb example of Carter's style of clarity and his way of referencing 
the axioms for the cohomology spaces $H_c^i(X)$.

Carter's book \cite{MR1266626} covers the character theory of finite 
Chevalley groups up until about 1984, including a quasi-exhaustive 
bibliography on the subject. In Chapter~12, the main results of Lusztig's 
book \cite{MR0742472} (which appeared in 1984) are only surveyed; these 
give explicit formulae for the decomposition of the generalised characters 
$R_{T,\theta}$ as linear combinations of irreducible characters, thereby
providing a complete classification of the irreducible characters. A further,
extremely helpful feature of \cite{MR1266626} is Chapter~13, which collects
a variety of explicit information, partly in the form of tables,  across 
all the topics touched upon in the book (e.g., tables of unipotent classes, 
degrees of unipotent characters and so on). Since the appearance of
\cite{MR1266626}, the whole subject has seen still further substantial new 
developments. For more recent accounts see the books by Bonnaf\'e 
\cite{MR2274998}, Cabanes--Enguehard \cite{MR2057756}, Digne--Michel 
\cite{MR4211777} and Geck--Malle \cite{MR4211779}.

\section{Kazhdan--Lusztig polynomials and cells} \label{sec-kl}

In 1979, Kazhdan and Lusztig published their fundamental paper 
\cite{MR0560412}, by which the authors were able to explain
many combinatorial phenomena that are associated with a finite 
Weyl group or, more generally, an arbitrary Coxeter group~$W$. The
key discovery is that of certain polynomials $P_{y,w}$ ($y,w\in W$) in
one variable with integral coefficients. These are used to:
\begin{itemize}
\item define partitions of $W$ into left, right and two-sided cells;
\item construct certain representations of the (generic) Iwahori--Hecke
algebra $H$ of $W$;
\item give a conjectural formula for certain (hitherto unknown) multiplicities 
in the Jordan--H\"older series of Verma modules (of a Lie 
algebra $\mathfrak{g}$ with Weyl group~$W$).
\end{itemize}
This has grown into a huge theory, with a significant literature touching 
upon a wide spectrum of areas in Lie theory; see also Lusztig's own comments 
\cite{Lcomments} for further background. Roger Carter does not seem to 
have publications directly concerned with Kazhdan--Lusztig polynomials or
cells. But he certainly followed the development of that theory very closely;
he also gave lectures and seminar talks about it at various occasions, 
thereby contributing to the dissemination of those important ideas and
stimulating further research. We mention two concrete examples.

The first-named author (M.G.) can actually speak from his own experience,
since he spent the academic year 1988/89 at Warwick University as a 
post-graduate visiting student, with Roger as supervisor. During that year 
Roger gave a one-semester lecture course on groups with a $(B,N)$-pair,
which included a chapter on Hecke algebras, Kazhdan--Lusztig polynomials 
and cells. This was followed by a one-semester seminar in which he 
presented the main ideas and constructions of Lusztig's series of papers 
on cells in affine Weyl groups (e.g., the $a$-function and the 
``asymptotic algebra'' $J$ of~$W$). Roger also gave a seminar talk about
the above-mentioned conjecture and the completion of its proof by 
Beilinson--Bernstein and Brylinski--Kashiwara. 

\begin{displayquote} {\small
``I was very fortunate to learn~---~still a student~---~about that amazing 
and extremely powerful work of Kazhdan--Lusztig and Lusztig from Roger's 
lectures. Being already familiar with it was enormously helpful in some of 
my later work, e.g., in \cite{MR2336039} where the technology of Lusztig's 
papers is used to show certain structural properties of Hecke algebras.'' 
\qquad (M.G.)}
\end{displayquote}

As a second example, we have reports from Jian-Yi Shi who became
Roger's PhD student in 1981. According to Shi, the story began as 
follows:
\begin{displayquote} {\small 
``Carter asked me to read the paper of Kazhdan and Lusztig \cite{MR0560412}, 
in October 1981. In order to test if I could well understand the main 
concept of the paper, he asked me to describe all the Kazhdan--Lusztig 
cells explicitly for the affine Weyl group $W_a(\tilde{A}_2)$. At that 
stage, he wouldn’t tell me that Lusztig had already described the cells
for all the affine Weyl groups of rank~$2$. In December 1981, I was able 
to describe explicitly all the cells of the affine Weyl group 
$W_a(\tilde{A}_2)$.''}
\end{displayquote}
Shi's thesis was completed in 1984. The main result is a precise description 
for all the left/right cells of the affine Weyl groups $W_a(\tilde{A}_n)$ 
(any $n\geq 2$), and also of the two-sided cells using Lusztig's $a$-function.
This was a big breakthrough in the theory of cells for affine Weyl groups
of type $\tilde{A}_n$. As a further outcome, Shi was led to define the 
concept of ``sign types'' for any affine Weyl group $W_a$, which is related 
to the partition of $W_a$ into cells. He also found an explicit formula 
for the number of ``sign types'' for $W_a(\tilde{A}_n)$. Then Roger  
formulated a conjecture about the number of ``sign types'' for arbitrary 
$W_a$; this conjecture was subsequently proved by Shi \cite{MR0871765}
(where Roger is explicitly credited for the formulation of the conjecture).

The results of Shi's thesis itself were considered to be suitable for
publication in the Springer LMS series. However, it was suggested that he 
provide more detailed background on the topics covered in the thesis. Shi
remembers:
\begin{displayquote} {\small 
``I revised my PhD thesis by adding three chapters under the guidance
 of Carter according to the referee’s comments. Meanwhile, Carter 
kindly did all the corrections for types, he also drew many graphs and
figures by his hands which couldn’t be properly typed at that time. This was 
a time-consuming hard work which needed much care and patience. Here I would 
like to express my deep gratitude for his generous help.''}
\end{displayquote}
Shi's work was published in 1986, see \cite{MR0835214}.

\section{Quantised enveloping algebras and canonical bases} 
\label{sec-further}

Throughout his career, both in his research and his teaching, Roger 
Carter had a keen interest in the classical Cartan--Killing theory of 
simple Lie algebras over $\CC$, their representations and related 
subjects. As Ian Stewart (his colleague at Warwick) relates:
\begin{displayquote}
{\small ``Soon after his appointment to the University of Warwick, Roger 
produced an elegant set of (unpublished) lecture notes on the classical
theory of structure and representations of finite-dimensional Lie algebras 
over algebraically closed fields of characteristic zero. These were issued 
in stapled, duplicated form with the iconic Warwick Mathematics Department’s 
orange covers.\footnote{Much later this material would be integrated in 
Carter's monograph \cite{MR2188930}.} When I taught an MSc course on 
Lie algebras, I drew heavily on Roger’s notes. Like all of his work, they 
were beautifully constructed, organised, and concise. They saved me a lot 
of work.''}
\end{displayquote}
So again let $\fg$ be a finite-dimensional simple
Lie algebra over $\CC$. We now consider its enveloping algebra $U(\fg)$; 
this is an infinite-dimensional associative algebra which plays an important 
role for the representation theory of~$\fg$. In the 1980s a new development
in this classical area emerged when Drinfel'd and Jimbo independently 
introduced a remarkable Hopf algebra deformation $U_v(\fg)$ of $U(\fg)$, 
which depends on a formal parameter~$v$ and subsequently became known as 
``quantum group'' or ``quantised enveloping algebra''. Roger had spent 1991 
in the United States visiting MIT to work with G.~Lusztig, who had just 
discovered the ``canonical basis'' of $U_v(\fg)$ (or, more precisely, of 
the negative part of $U_v(\fg)$ in a triangular decomposition);  for more
details~---~and also relations to the geometry of quiver varieties and to
Kashiwara's work on similar bases~---~see Lusztig's widely cited book 
\cite{MR2759715} on the subject. The ``canonical basis'' is hard to compute 
but a truly marvelous object. One of its features is that, by specialising 
the formal parameter $v$ to~$1$, it simultaneously gives rise to ``canonical 
bases'' in all finite-dimensional irreducible representations of $\fg$
itself. Thus, although the existence of the canonical basis could only be 
seen in the context of $U_v(\fg)$ or, alternatively, via the geometry of 
quiver varieties, this basis has a striking impact on the classical
representation theory of the original Lie algebra~$\fg$. For example, the
adjoint representation of $\fg$ (that is, $\fg$ considered as a module over 
itself with operation given by the Lie bracket) certainly is a 
finite-dimensional irreducible representation. The canonical basis in this 
case is explicitly described by Lusztig \cite{MR3589909}. (We already 
encountered this reference before in connection with ``canonical'' Chevalley 
bases for~$\fg$; see Section~\ref{sec-two}, p.~\pageref{candj}.) 

It will come as no surprise that Roger immediately became interested in
these exciting new developments. Back at Warwick, he gave four new PhD 
students, J.~Cockerton, R.~Green, B.~Marsh and A.~Scott, projects to
work on relating to canonical bases.

B.~Marsh remembers:
\begin{displayquote}{\small ``Keeping up to date with the latest 
developments was one of Roger’s strengths. My thesis problem was 
to compute the canonical basis of a quantum group explicitly. There 
was plenty of scope for partial results, since the canonical basis induces 
bases (also referred to as canonical bases) in all finite-dimensional 
irreducible representations. I was able to find canonical bases in some 
modules and find a relationship with standard monomial theory. A 
particularly nice outcome of Roger’s approach to supervision was a joint 
paper with Richard Green \cite{MR1604827}. After my thesis, Roger 
and I stayed in touch and we wrote joint papers \cite{MR1800743}, 
\cite{MR1903958} on the combinatorics of parametrizations of the 
canonical basis arising from the approaches of Kashiwara and Lusztig.''}
\end{displayquote}

In addition to the joint publications with Marsh, Roger wrote a single 
author paper \cite{MR1635675} exploring Lusztig's piecewise linear function 
for computing the canonical basis and applies the algorithm to certain 
small rank Lie algebras. He also wrote a survey article \cite{MR1670769} 
which gives an introduction to all these striking ideas, as well as to 
Kashiwara's work on ``crystal bases'' and Littelmann's ``path model''
approach. 

B.~Marsh goes on to say:
\begin{displayquote}
{\small ``I will never forget the many good times I enjoyed being
supervised by Roger, and collaborating and interacting with him. Roger was 
always laughing and enjoying mathematics, and I have learnt a lot from 
working with him. He has had a profound impact on my career, and on Lie 
theory more generally.''}
\end{displayquote}
And as A.~Scott remembers:
\begin{displayquote}
{\small ``Roger regularly took his four 1992 intake of PhD students 
for pub lunches and always picked up the bill. Directors bitter was his 
preferred drink on these occasions, though recollections may vary. 

I still have the beautiful letter Roger wrote to me (in his perfectly
upright handwriting, much like the man himself), congratulating me upon 
submission of my thesis, whilst oh-so-gently disabusing me of any notion
that it was worth circulating more widely!''}
\end{displayquote}

\section{Concluding remarks} \label{sec-conclude}

Throughout we have tried to focus on subjects that~---~to the best of our 
knowledge~---~have been long-running themes in Roger Carter's life as a 
mathematician: soluble groups in the early stages of his career, and then 
Lie algebras, algebraic groups, representation theory and Lie theory in 
general. There are a few topics that we omitted from the discussion, partly
because we did not feel sufficiently competent to handle them in detail. 
As an example, we mention the theory of Kac--Moody groups, which are
infinite-dimensional generalisations of the Chevalley groups that we 
encountered earlier in this text. Roger had a PhD student, Jacqui Ramagge, 
working on that subject. Y. Chen relates:

\begin{displayquote} 
{\small ``From 1988 to early 1990, I worked in Warwick as a Research
Fellow under the supervision of Roger, who became a dear friend of me.''}
\end{displayquote}
Roger encouraged Chen to work on quantum groups and/or Kac--Moody theory.
Eventually, their collaboration resulted in the substantial joint paper 
\cite{MR1206622} about the automorphism groups of Kac--Moody groups, in 
which they proved a conjecture of Kac and Peterson from the 1980s for a 
large class of those groups. This was later extended to all possible types 
by Caprace--M\"uhlherr \cite{MR2180452}. P.-E. Caprace remembers a letter 
that he received from Roger after sending him a draft of \cite{MR2180452}:

\begin{displayquote} 
{\small ``I was very suprised at the time that he took the trouble to 
write this letter to express interest and encouragement, although I did 
not know him personally. Retrospectively I realise that it witnessed a 
remarkable generosity and a genuine care for other members of the 
mathematical community.''}
\end{displayquote}

In addition to his articles, books and other written work, Roger had a 
great influence on many people through the numerous seminar talks and 
lecture courses that he gave at various occasions. Burkhard K\"ulshammer 
took notes of such a lecture course given by Roger in 1987. He remembers:

\begin{displayquote} {\small ``On January 22 he started his series of 
lectures in Essen. He lectured once each week, for 90 minutes, until 
April 1. I was then a postdoc and ``Privatdozent'' in Dortmund and came 
to Essen to attend his lectures and to take notes. After each lecture 
I prepared the notes electronically, and in the following week Roger 
corrected them. My task was relatively easy since his lectures 
were~---~as always~---~crystal clear and could usually be copied word 
by word.  The title of his lecture series was "On the structure and
complex representation theory of finite groups of Lie
type", and the contents was mainly an overview of his
corresponding book(s).''}
\end{displayquote}

In order to make the point about Roger Carter's continuing, wide interests 
and his profound expertise in various areas linked to algebraic groups and 
Lie theory in general, we conclude by only mentioning here some publications 
that he wrote later in his career: there is the substantial monograph on Lie
algebras of finite and affine type \cite{MR2188930} (published in 2005), 
his superb survey on the work of G. Lusztig \cite{MR2235338} (2006), and 
also his introductions to ``monstrous moonshine'' \cite{MR1951481} (2002)
and the relatively recent Fomin--Zelevinsky theory of ``cluster algebras''
\cite{MR2267299} (2006). About the latter, Ana Paula Santana remembers:
\begin{displayquote}
{\small ``I first met Roger Carter in 1986, when I was a PhD student in 
Warwick. I remember well his characteristic laughter and what an amazing 
mathematics communicator he was. Later, after I came to the University of 
Coimbra, Roger visited our Department and gave two excellent advanced 
courses, on ``Quantised enveloping algebras'', in 1995, and 
``Cluster Algebras'', in 2003.''}
\end{displayquote}
All this shows his genuine desire and a unusual ability to share
his knowledge, his expertise~---~and his joy and excitement about new 
developments~---~with his students, his colleagues, and with the general 
mathematical community.

\medskip
{\bf Acknowledgements}. We are indebted to a number of colleagues for
sending us comments about Roger's work, including remarks on history,
as well as impact on and connections to subsequent work. We could not
incorporate all of these comments in this text, but in any case they
helped us to put things into perspective. Special thanks go to
Pierre--Emmanuel Caprace, Yu Chen, Michael Collins, Steve Donkin, Paul Fong,
Trevor Hawkes, Gerhard Hiss, Bob Howlett, Sasha Kleshchev, Burkhard
K\"ulshammer, George Lusztig, Richard Lyons, Gunter Malle, Bethany Marsh,
Andrew Mathas, Dmitriy Rumynin, Jian-Yi Shi.

\bigbreak

\bibliographystyle{amsplain} 

\providecommand{\MR}{\relax\ifhmode\unskip\space\fi MR }
\providecommand{\MRhref}[2]{%
  \href{http://www.ams.org/mathscinet-getitem?mr=#1}{#2}
}
\providecommand{\href}[2]{#2}

\end{document}